\documentstyle[11pt]{article}

\title{Pavelka-style completeness in expansions of \L ukasiewicz logic}

\author{{\sc Hector Freytes} \thanks{Permanent address: Instituto Argentino de Matem\'atica (IAM), Saavedra 15 - 3er Piso - 1083  Buenos Aires, Argentina, CONICET.}  \thanks{The author express his gratitude to Roberto Cignoli, for his advice during the preparation of this paper.}}

\date{{\small University of Cagliari, \\
Dipartimento di Scienze Pedagogiche e Filosofiche, \\ Via Is Mirrionis 1,09123,
Cagliari-Italia \\e-mail: hfreytes@gmail.com }}

\begin{document}

\bibliographystyle{plain}

\maketitle

\begin{abstract}

\noindent An algebraic setting for the validity of Pavelka style
completeness for some natural expansions of \L ukasiewicz logic by
new connectives and rational constants is given. This algebraic
approach is based on the fact that the standard MV-algebra on the
real segment $[0, 1]$ is an injective MV-algebra. In particular the
 logics associated with MV-algebras with product and with
divisible MV-algebras are considered. \vspace{0,3cm}


\end{abstract}

\begin{small}

\end{small}

\bibliography{pom}

\newtheorem{theo}{Theorem}[section]

\newtheorem{definition}[theo]{Definition}

\newtheorem{lem}[theo]{Lemma}

\newtheorem{prop}[theo]{Proposition}

\newtheorem{coro}[theo]{Corollary}

\newtheorem{exam}[theo]{Example}

\newtheorem{rema}[theo]{Remark}{\hspace*{4mm}}

\newtheorem{example}[theo]{Example}

\newcommand{\proof}{\noindent {\em Proof:\/}{\hspace*{4mm}}}

\newcommand{\qed}{\hfill$\Box$}

\section*{Introduction}
The completeness theorem for fuzzy propositional logic was proved by
J. Pavelka in \cite{PAV1}, who built a propositional many valued
logical systems over \L ukasiewicz logic adding to the language a
truth constant $\overline r$ for each $r \in [0,1]$, together with
 additional axioms. Later, H\'ajek gave a proof of the Pavelka
completeness of \L ukasiewicz logic with rational constants, based
on the continuity of \L ukasiewicz implication on the real segment
$[0, 1]$ (see \cite[\S3.3]{HAJ}). Following the ideas of  H\'ajek's
proof, analogous results were obtained for some expansions of \L
ukasiewicz logic (see \cite{EGM}, \cite{GER} and \cite{HC}).

The aim of this paper is to offer a new method for proving the
validity of Pavelka style completeness of \L ukasiewicz logic with
rational constants, based on the crucial fact that the standard
MV-algebra $[0, 1]_{MV}$ is injective in the category of
MV-algebras. This method can also be applied for a special class of
expansions of \L ukasiewicz logic, called compatible expansions.
Such expansions guarantee the injectivity of the standard MV-algebra
$[0, 1]_{MV}$ with the corresponding additional operations.

The paper is organized as follows. In section 1 we recall some basic
definitions and properties of injectives, algebraizable logics and
compatible expansions.  Section 2 contains the main results of the
paper. It is shown  that Pavelka completeness can be established
under some mild hypothesis on compatible expansions of \L ukasiewicz
logic. A type compactness is also proved. The results of section 2
are applied in section 3 to obtain a new proof of Pavelka
completeness for \L ukasiewicz logic, as well as  Pavelka style
completeness for product \L ukasiewicz logic and for divisible \L
ukasiewicz logic.

\section{Basic notions}
We use freely all basic notions of universal algebra that came be
found in \cite{Bur}. Let ${\cal A}$ be a class of algebras of type
$\tau$. For all algebras $A, B$ in ${\cal A}$, $[A,B]_{\cal A}$ will
denote the set of all homomorphism $g:A\rightarrow B$.

An algebra $A$ in ${\cal A}$ is {\it injective} iff for every
monomorphism  $f \in [B, C]_{{\cal A}}$ and every  $g \in
[B,A]_{{\cal A}}$ there  exists  $h \in [C, A]_{{\cal A}}$ such that
$hf = g$. $A$ is {\it self-injective} iff every homomorphism from a
subalgebra of $A$ into $A$, extends to an endomorphism of $A$.

Recall from \cite{FRE} that a  simple algebra $I_M$ is said to be
{\it maximum simple} iff for each simple algebra $I$, $I$ can be
embedded into $I_M$. Two constant terms  $0,1$ of the language of
$\cal A$ are called  {\it distinguished constants} iff $A\models
0\not=1$ for each nontrivial algebra $A$ in $\cal A$.

\begin{theo}\label{Simple injective}
{\rm \cite[Theorem 3.4]{FRE}} Let ${\cal A}$ be a variety satisfying
the congruence extension property, with distinguished constants $0,
1$. If $I$ is a self-injective maximum simple algebra in ${\cal A}$
then $I$ is injective. \qed
\end{theo}

If $\tau$ is a family of finitary function symbols, by a {\it logic}
${\cal L}$ of type  $\tau$ we will understand a structural finitary
consequence relation $\vdash_{\cal L}$ in the absolutely free
algebra $Fm_{\tau}$ of type  $\tau$ generated by the proposotional
variables $p_1, p_2, \ldots$. Each function symbol of $\tau$ is
called a {\it connective}, and the terms of  $Fm_{\tau}$ are called
{\it formulas}. Usually, $\vdash_{\cal L}$ is specified by a set of
Hilbert style axiom schemes and inference rules. A set $T$ of
formulas is called a {\it theory} of ${\cal L}$ if  $T \vdash_{\cal
L} \varphi$ implies $\varphi \in T$ for every $\varphi \in
Fm_{\tau}$.\footnote{Notice that our notion of a theory differs from
the one used in \cite{HAJ}, where a theory is just a set of
formulas.} We denote by $Th{\cal L}$ the lattice of theories
associated to ${\cal L}$ {\rm (\cite[1.1]{BlokPig1})}. If $S$ is a
set of formulas, we denote by $T(S)$ the {\it  theory generated} by
$S$. Recall {\rm \cite[Definition 2.8, Corollary 2.9]{BlokPig1}}
that a logic ${\cal L}$ is {\it algebraizable} by a quasivariety $K$
of type $\tau$  via a finite set of binary formulas $p
\Leftrightarrow q =(p \Leftrightarrow_i q) _i$ (called the
equivalence formulas of ${\cal L}$) and a finite set of identities $
\delta (p) \approx \epsilon (q) = (\delta_j(p) \approx
\epsilon_j(q))_j$ (the defining equation of ${\cal L}$) if and only
if the following conditions are met:

\begin{enumerate}
\item[i]
$\phi_1, \ldots, \phi_n \vdash_{\cal L} \phi$ iff $\{\delta(\phi_i)
\approx \epsilon(\phi_i)\}_{i=1}^n \models_K \delta(\phi) \approx
\epsilon(\phi) $

\item[ii]
$\varphi \approx \psi  \models_K \delta (\varphi \Leftrightarrow
\psi) \approx \epsilon (\varphi \Leftrightarrow \psi)$

\end{enumerate}

\noindent we are using the abbreviations $p \Leftrightarrow q$ for
$\{p \Leftrightarrow_i q \}_i$, $\delta (p \Leftrightarrow q)
\approx \epsilon (p \Leftrightarrow q) $ for $ \{ \delta_j (p
\Leftrightarrow_i q) \approx \epsilon_j (p \Leftrightarrow_i q)
\}_{i,j}$. The quasivariety ${\cal K}$, is uniquely determined by
${\cal L}$, when it exists and the equivalence formulas and defining
equations are also unique in the sense that for any other system $p
\Leftrightarrow q$, $\delta'(p) \approx \epsilon' (q)$ satisfing
{\it i, ii}, $\vdash p \Leftrightarrow q $ iff $\vdash p
\Leftrightarrow q $ and $\models_K \delta'(p) \approx \epsilon' (q)$
iff $\models_K \delta(p) \approx \epsilon(q)$. Will be denote this
quasivariety as  ${\cal K}_{\cal L}$. From {\rm \cite[Theorem 4.7
and Theorem. 4.10 ]{BlokPig1}}, if we consider for each $T \in
Th{\cal L}$

$$\Omega T = \{ \varphi \approx \psi : \varphi \Leftrightarrow \psi \in T \} $$

\noindent then $\Omega T$ is a congruence in $Fm_\tau$ and the
quasivariety  generated by the set $\{Fm_\tau/\Omega T: T \in Th
{\cal L} \}$ is ${\cal K}_{\cal L}$. In what follows, for each
$\alpha \in Fm_\tau$ will be denoted by $[\alpha]$ the equivalence
class of $\alpha$ with respect $Fm_\tau/\Omega T$. Let $A\in {\cal
K}$, a valuation of formulas on $A$ is an homomorphism $v: Fm_{\tau}
\rightarrow A$. If $T \in Th{\cal L}$ and $a\in A$, $v(T) = a$
means, as usual, that  $v(\alpha )= a$ for each $\alpha \in T$.

Let $E$ be a set of equations in the vocabulary $\tau$ that gives an
axiomatization of a variety ${\cal A}$. If $\sigma = (f_i)_{i \in
I}$ is a family of operation symbols such that $\sigma \cap \tau =
\emptyset$ and $E(\sigma)$ a set of equation in the expanded
vocabulary $\tau \cup \sigma$, we denote by ${\cal A}_{E(\sigma )}$
the variety of type $\tau \cup \sigma$ defined by the equations in
the set $E \cup E(\sigma)$. We will refer to ${\cal A}_{E(\sigma )}$
as the $E(\sigma )$-expansion of ${\cal A}$. For each $A \in {\cal
A}_{E(\sigma )}$, we denote by $Con_{\cal A} (A)$ the lattice of
${\cal A}$-congruences of $A$ and by $Con_{{\cal A}_{E(\sigma
)}}(A)$ the lattice of ${\cal A}_{E(\sigma )}$-congruences of $A$.
We say that ${\cal A}_{E(\sigma)}$ is a {\it compatible expansion}
of ${\cal A}$ iff for each $A \in {\cal A}_E{\sigma}$ $Con_{{\cal
A}_{E(\sigma )}}(A) = Con_{\cal A} (A)$. Note that  an algebra is
simple in the compatible expansion ${\cal A}_{E(\sigma )}$ iff  it
is simple in ${\cal A}$.

Let $A \in {\cal A}$, then $A$ admits  an $E(\sigma)$-expansion iff
there is a family of $A$-operation $\sigma_{A} = (f_i^A)_{i\in I}$
making $(A, \sigma_A)$ into a ${\cal A}_\sigma$-algebra. This
$E(\sigma)$-expansion on $A$ is {\it canonical} iff $\sigma_{A}$ is
unique and each sub algebra of $A$ admits at most one
$E(\sigma)$-expansion.

\begin{prop}\label{Simple EXP}
Let ${\cal A}$ be a variety and ${\cal A}_{E(\sigma)}$ be a
compatible expansion. If $I$ is maximum simple algebra in ${\cal A}$
admiting a canonical $E(\sigma)$-expansion, then $I$ is maximum
simple in ${\cal A}_{E(\sigma)}$.
\end{prop}

\begin{proof}
Let $(I , \sigma_I)$ be canonical  $E(\sigma)$-expansion of $I$ and
$A$ be a simple algebra in ${\cal A}_{E(\sigma)}$. Since $A$ is
simple in ${\cal A}$, there exists a sub ${\cal A}$-algebra $A_0$ of
$I$ ${\cal A}$-isomorphic to $A$. An ${\cal A}$-isomorphism
$i:A\rightarrow A_0$  induces $\sigma$-operations on $A_0$ making
$(A_0 , \sigma_{A_0})$ into a ${\cal A}_{E(\sigma)}$-algebra. Since
the expansion $(I , \sigma_I)$ is canonical,  $\sigma_{A_0}$ is the
set of  restrictions to $A_0$ of the $\sigma_I$-operations . Thus
$i$ is a ${\cal A}_\sigma$-monomorphism, proving that $I$ is maximum
simple algebra in ${\cal A}_\sigma$.

\qed
\end{proof}

\section{MV-expansions and Pavelka-style completeness}
\begin{definition}
An MV-algebra {\rm \cite{CDM}, \cite{HAJ}} is an algebra $ \langle
A, \land, \lor, \odot, \rightarrow, 0, 1 \rangle$ of type $ \langle
2, 2, 2, 2, 0, 0 \rangle$ satisfying the following axioms:

\begin{enumerate}
\item
$\langle A,\odot,1 \rangle$ is an abelian monoid,

\item
$L(A) = \langle A, \lor, \land, 0,1 \rangle$ is a bounded lattice,

\item
$(x \odot y)\rightarrow z = x\rightarrow (y\rightarrow z)$,

\item
$((x\rightarrow y)\odot x)\land y = (x\rightarrow y)\odot x$,

\item
$(x\land y)\rightarrow y = 1$,

\item
$x\odot (x\rightarrow y) = x\land y$,

\item
$(x\rightarrow y)\lor (y\rightarrow x)  = 1$,

\item
$(x\rightarrow 0)\rightarrow 0 = x$

\end{enumerate}
\end{definition}

\noindent We denote by ${\cal MV}$ the variety of MV-algebras. In
agreement with the usual MV-algebraic operations we define the {\it
negation} as the unary operation $\neg x = x\rightarrow 0$. and the
binary operation $x \oplus y = \neg (\neg x \odot  \neg y)$,
$x\rightarrow y = \neg x \oplus y$. For element $x$ in a MV-algebra
and $n\in N$, we denote $nx$ the element inductively defined by $0x
= 0$, $(n+1)x = (nx) \oplus x$. ${\cal MV}$ satisfies congruence
extension property {\rm  \cite[Theorem 1.8]{BlokFer}} and then each
compatible expansion also satisfies. An important example is
$[0,1]_{MV}= \langle [0,1], \odot, \rightarrow, \land, \lor, 0, 1
\rangle$ such that $[0,1]$ is the real unit segment,$\land$, $\lor$
are the natural meet and join on $[0,1]$ and $\odot$ and
$\rightarrow$ are defined as follows: $x\odot y:= max(0,x+y-1)$,
$x\rightarrow y:= min(1,1-x+y)$. $[0,1]_{MV}$ is the maximum simple
algebra in ${\cal MV}$ {\rm (see \cite[Theorem 3.5.1]{CDM})}.
Moreover $[0,1]_{MV}$ is a rigid algebra {\rm (see \cite[Corollary
7.2.6]{CDM})}, hence self-injective resulting injective in the
variety ${\cal MV}$ in view of Theorem \ref{Simple injective}.

\begin{prop} \label{INJEXP}
Let ${\cal MV}_{E(\sigma)}$ be a compatible expansion of ${\cal MV}$
such that $[0,1]_{MV}$ admits a canonical $E(\sigma)$-expansion.
Then the $E(\sigma)$-expansion of $[0,1]_{MV}$ is injective in ${\cal
MV}_{E(\sigma)}$.
\end{prop}

\begin{proof}
Let $[0,1]_{E(\sigma)}$ be the $E(\sigma)$-expansion of
$[0,1]_{MV}$. Since ${\cal MV}_{E(\sigma)}$ is a compatible
expansion of ${\cal MV}$, it satisfies the congruence extension
property, and clearly $0,1$ are distinguished constant. Then by
Proposition \ref{Simple EXP} $[0,1]_{E(\sigma)}$ is a maximum simple
algebra. Moreover $[0,1]_{E(\sigma)}$ is self injective since the
expansion is canonical. Thus by Theorem  \ref{Simple injective}
$[0,1]_{E(\sigma)}$ is injective.

\qed
\end{proof}

\begin{definition}
{\rm
 A compatible expansion ${\cal MV}_{E(\sigma)}$ of ${\cal MV}$
is said to be {\it admissible} if

\begin{enumerate}
\item
$[0,1]_{MV}$ admits canonical $E(\sigma)$-expansion,

\item
$Q_{[0,1]}$ can be $E(\sigma)$-expanded, were $Q_{[0,1]}$ denotes
the subalgebra of rational numbers in $[0,1]_{MV}$.

\end{enumerate}
}
\end{definition}

\begin{definition}
{\rm We can define the \L ukasiewicz propositional calculus $\L$
{\rm \cite{CDM}} from $\tau = \{\rightarrow, 0 \}$ where
$\rightarrow$ is a binary connective and $0$ is constant. Futher
connectives are defined as follows: $\neg \alpha$ is $\alpha
\rightarrow 0$, $1$ is $\neg 0$ and $\alpha \odot \beta$ is $\neg
(\alpha \rightarrow \neg \beta)$. The following formulas are axioms

\begin{enumerate}

\item[\L 1]
$\alpha \rightarrow (\beta \rightarrow \alpha)$

\item[\L 2]
$(\alpha \rightarrow \beta) \rightarrow ((\beta \rightarrow \gamma )
\rightarrow (\alpha \rightarrow \gamma)) $

\item[\L 3]
$(\neg \alpha \rightarrow  \neg \beta) \rightarrow (\beta
\rightarrow \alpha) $

\item[\L 4]
$((\alpha \rightarrow \beta) \rightarrow \beta) \rightarrow ((\beta
\rightarrow \alpha) \rightarrow \alpha)$

\end{enumerate}

\noindent The unique deduction rule is the {\it Modus Ponens}

}
\end{definition}

This calculus is  algebraizable in the variety $\cal MV$ via the
system $p \Leftrightarrow q = \{p\rightarrow q, q\rightarrow p \}$
and $\delta(p) \approx \epsilon(p) = \{p \approx p\rightarrow p \} $
(see \cite{RTV}).

\begin{definition}
{\rm By an {\it admissible expansion} of \L ukasiewicz logic $\L$ we
understand a logic $\L_{\sigma}$ involving a set $\sigma$ of new
connectives and new axioms, algebraizable in an admissible expansion
${\cal MV}_{E(\sigma)}$ of ${\cal MV}$. }
\end{definition}

\begin{prop} \label{COMPLETE}
Given an  admissible expansion $\L'$  of $\L$, an $\L'$-formula
$\alpha$, and $T \in Th_{\L'}$ such that $\alpha \not \in T$, there
exists $T' \in Th_{\L'}$ such that $T \subseteq T'$, $\alpha \not
\in T'$ and $Fm_{\L '}/\Omega T'$ is totaly ordered.
\end{prop}

\begin{proof}
Follows from {\rm \cite[Lemma 2.4.2]{HAJ}} and the fact that
compatible expansions do not involve new inference rules.
\qed \\
\end{proof}

\begin{definition}{\rm\cite[\S3.3]{HAJ}}
{\rm A $\sigma$- {\it Pavelka expansion } of \L ukasiewicz logic is
obtained by adding into the language of an admissible expansion
$\L_{\sigma} $ of $\L$,  {\it truth constant} $\overline r$ for each
$r\in Q\cap [0,1]$, together with the following additional {\it
book-keeping} axioms:

\begin{enumerate}

\item
$\overline 0 \Longleftrightarrow  0$,

\item
$ \overline r \rightarrow \overline s \Longleftrightarrow \overline
{r \rightarrow s}$,

\item
$f(\overline r_1 \ldots \overline r_n ) \Longleftrightarrow
\overline {f(r_1 \ldots r_n )}$ for each $n$-ary connective $f$ in
$\sigma$ and $r_1 \ldots r_n \in Q$.

\end{enumerate}

\noindent Moreover, for each theory  $T$  and  each formula $\alpha$
 in the $\sigma$- Pavelka expansion we define:

\begin{enumerate}

\item[]
The {\it truth degree} of $\alpha$ over T is $\Vert \alpha \Vert_T =
\bigwedge \{v(\alpha):v(T)= 1 \mbox{ {\it in} $[0,1]$} \}$.

\item[]
The {\it proof degree} of $\alpha$ is $\vert \alpha \vert_T =
\bigvee \{ r \in Q: \overline r \rightarrow \alpha \in T \}$.
\end{enumerate}
}
\end{definition}

Observe that $\sigma$-Pavelka expansion is an admissible expansion
$\L_{\sigma \cup \{\overline {r} : \, r \in Q\}} $ of $\L$.

We denote by $Fm_{\sigma}( \overline Q)$ the set of the formulas of
the $\sigma$-{\it Pavelka expansion}. The following theorem is a
generalization of  Pavelka completeness for \L ukasiewicz logic (see
\cite[\S3.3]{HAJ}):

\begin{theo}\label{COMP}
Let $\L_{\sigma}$ be a admissible expansion of $\L$. If $T$ is a
theory and $\varphi$ is a formula, both in the $\sigma$- Pavelka
expansion of $\L$, then

$$\vert \varphi \vert_T   = \Vert \varphi \Vert_T.$$
\end{theo}

\begin{proof}
Let $[0,1]_{E(\sigma)}$ the $E(\sigma)$-expansion of $[0,1]_{MV}$
and the $Q_{E(\sigma)}$ the $E(\sigma)$-expansion of $Q_{[0,1]}$.
Assume that $T \neq Fm_{\sigma} (\overline Q) $. We first prove that
$\vert \varphi \vert_T$ is a lower bound of $\{v(\varphi):v(T)= 1
\}$. Let $v$ be a valuation such that $v(T)=1$. If $r \in  Q$ is
such that $\overline r \rightarrow \varphi \in T$ then, we have that
$1 = v(\overline r \rightarrow \varphi) = r \rightarrow v(\varphi)$.
Thus $r\leq v(\varphi)$ resulting $\vert \varphi \vert_T \leq
v(\varphi)$. We proceed now to prove that $\vert \varphi \vert_T$ is
the greatest lower bound of $\{v(\varphi):v(T)= 1 \}$. In fact, let
$b$ be a lower bound of $\{v(\varphi):v(T)= 1 \}$. Suppose that
$\vert \varphi \vert_T < b$. Then there exists $r_0 \in Q$ such that
 $\vert \varphi \vert_T < r_0 < b$. Thus $T$ does not prove $\overline
{r_0} \rightarrow \varphi$. By Proposition \ref{COMPLETE} there
exists a  theory $T'$ such that $T \subseteq T'$, $T'$ does not
prove $\overline  {r_0} \rightarrow \varphi$ and
$Fm_{\sigma}(\overline Q) /\Omega T'$ is a totaly ordered algebra in
${\cal MV}_{E(\sigma)}$  containing  $\{[\overline r]\}_{r \in Q}$
as sub algebra isomorphic to $Q_{E(\sigma)}$. Since $T'$ does not
prove $\overline {r_0} \rightarrow \varphi$, $[\overline r_0]
\rightarrow [\varphi] < 1$ in $Fm_{\sigma}(\overline Q) /\Omega T'$
resulting $ [\varphi] < [\overline r_0]$ since it is a totaly order.
Let $i_1$ the canonical embeding $Q_{E(\sigma)} \rightarrow [0,1]_{E(\sigma)}$
and $i_2$ the
canonical embedding $Q_{E(\sigma)} \rightarrow Fm_{\sigma}(\overline Q) /\Omega
T'$. Since  $[0,1]_{E(\sigma)}$ is injective in the ${\cal MV}_{E(\sigma)}$,
there
exist an homomorphism $f:Fm_{\sigma}(\overline Q) /\Omega T'
\rightarrow [0,1]_{E(\sigma)}$ such that the following diagram is commutative

 \begin{center}
 \unitlength=1mm
 \begin{picture}(20,20)(0,0)
 \put(8,16){\vector(3,0){5}} \put(2,10){\vector(0,-2){5}}
 \put(10,4){\vector(1,1){7}}

 \put(2,10){\makebox(13,0){$\equiv$}}

 \put(2,16){\makebox(0,0){$Q_{E(\sigma)}$}}
 \put(23,16){\makebox(0,0){$[0,1]_{E(\sigma)}$}}
 \put(2,0){\makebox(0,0){$Fm_{\sigma}(\overline Q) /\Omega
 T'$}}
 \put(2,20){\makebox(17,0){$i_1$}} \put(2,8){\makebox(-6,0){$i_2$}}
 \put(18,2){\makebox(-4,3){$f$}}
 \end{picture}
 \end{center}

 By the commutativity $f([\varphi]) \leq r_0 < b$. If we consider the
 valuation $\pi:Fm_{\sigma}(\overline Q)  \rightarrow Fm_{\sigma}(\overline
 Q) /\Omega T'$ such that $\alpha \rightarrow [\alpha]$
 then the composition $fp$ is a valuation over $[0,1]_{E(\sigma)}$ such that
$fp(T') =
 1$ and $fp(T) = 1$ since $T\subseteq T'$ resulting $fp(\varphi)
 \in \{v(\varphi):v(T)= 1 \}$. But $fp(\varphi) \leq r_0 < b$ which is a
 contradiction since $b$ is a lower bound of $\{v(\varphi):v(T)= 1 \}$.
 Therefore $b \leq \vert \varphi \vert_T$ resulting $\vert \varphi
 \vert_T = \Vert \varphi \Vert_T$.
 \qed
 \end{proof}

\medskip

From the above completeness theorem, we can establish a kind of
compactness theorem.

\begin{theo} \label{COMPAC}
Let $S$ be a set of formulas and $\alpha$ be a formula in a
$\sigma$-{\it Pavelka expansion}.  Then we have:
$$\mbox{If} \hspace{0.2cm} r\leq \Vert \alpha \Vert_{T(S)} \hspace{0.2cm}
\mbox{then} \hspace{0.2cm} \exists \hspace{0.1cm} S_0 \subseteq S
\hspace{0.2cm} \mbox{finite} \hspace{0.1cm} \mbox{such that}
\hspace{0.2cm} r\leq \Vert \alpha \Vert_{T(S_0)}  $$

\end{theo}

\begin{proof}
If $r\leq \Vert \alpha \Vert_{T(S)}$ then by Theorem \ref{COMP},
$\overline r \rightarrow \alpha \in T(S)$. If  $\alpha_1, \cdots
\alpha_n, \overline r \rightarrow \alpha$ is a proof of $\overline r
\rightarrow \alpha$ from $S$, we can consider the finite set $S_0 =
\{ \alpha_1, \cdots \alpha_n \}$. Using again Theorem \ref{COMP} we
have $ r\leq  \vert \alpha \vert_{T_{S_0}}  =  \Vert \alpha
\Vert_{T_{S_0}}$ \qed
\end{proof}

\section{Applications}

\subsection{\bf \L ukasiewicz logic}
Since \L ukasiewicz logic $\L$, is an admissible expansion of itself
(with $\sigma = \emptyset$),   Theorem \ref{COMP} provides a new
proof of Pavelka completeness for \L ukasiewicz logic with rational
constants (cf. \cite[\S3.3]{HAJ}).

\subsection{\bf Product \L ukasiewicz logic}

\begin{definition}
{\rm A product MV-algebra {\rm \cite{MONT}} (for short: PMV-algebra)
is an algebra $ \langle A,  \bullet \rangle$ satisfying the
following

\begin{enumerate}
\item[1]
$A$ is an MV-algebra

\item[2]
$ \langle A, \bullet, 1 \rangle$ is an abelian monoid

\item[3] $x \bullet (y \odot \neg z) = (x \bullet y) \odot \neg (x \bullet z)
$

\end{enumerate}
}
\end{definition}

We denote by ${\cal PMV}$ the variety of PMV-algebras. If $A$ is a
PMV-algebra then {\rm  \cite[Lemma 2.11]{MONT}} $Con_{\cal PMV}(A) =
Con_{\cal MV}(A)$, resulting  ${\cal PMV}$ a compatible expansion of
${\cal MV}$. $[0,1]_{MV}$ equiped with the usual multiplication in
the unitary interval is a PMV-algebra denoted by $[0,1]_{PMV}$. It
contain $Q_[0,1]$ equiped with multiplication as sub PMV-algebra. It
is clear that $[0,1]_{PMV}$ is a hereditarily simple algebra. If $A$
is a semisimple MV-algebra then is at most one operation $\bullet$
making $ \langle A,  \bullet \rangle$ into a PMV-algebra {\rm
\cite[Lemma 3.1.14]{MR}}. Thus $[0,1]_{PMV}$ is the canonical
PMV-expansion of $[0,1]_{MV}$.

\begin{definition}
{\rm We define the product \L ukasiewicz propositional calculus
$P\L$ adding into the language of $\L$ the binary connective
$\bullet$ and considering the following formulas as axioms

\begin{enumerate}

\item[P\L 0]
\L ukasiewicz axioms

\item[P\L 1]
$ (\alpha \bullet \beta ) \rightarrow (\beta \bullet \alpha )$

\item[P\L 2]
$ (\top \bullet \alpha ) \Longleftrightarrow \alpha$

\item[P\L 3]
$ (\alpha \bullet \beta ) \rightarrow \beta$

\item[P\L 4]
$(\alpha \bullet \beta ) \bullet \gamma \Longleftrightarrow  \alpha
\bullet (\beta \bullet \gamma)$

\item[P\L 5]
$x \bullet (y \odot \neg z) \Longleftrightarrow  (x \bullet y) \odot
\neg (x \bullet z) $

\end{enumerate}

\noindent The unique deduction rule is the {\it Modus Ponens} }
\end{definition}

\begin{prop} \label{ALGPMV}
$P\L$ calculus is  algebraizable in ${\cal PMV}$.

\end{prop}

\begin{proof}
We only need to prove that $\alpha \Leftrightarrow \beta \vdash_{P
\L} \gamma \bullet \alpha  \Leftrightarrow \gamma \bullet \beta $.
For this we prove that $\vdash (\alpha \rightarrow \beta)\rightarrow
((\gamma\bullet \alpha) \rightarrow (\gamma \bullet \beta))$.

\begin{enumerate}

\item[(1)]
$\vdash \gamma \bullet (\alpha \odot \neg \beta) \rightarrow
((\alpha \odot \neg \beta))$ \hspace{3.2 cm} {\footnotesize by Ax
P\L 3}

\item[(2)]
$\vdash ((\gamma\bullet \alpha) \odot \neg(\gamma \bullet
\beta))\rightarrow \gamma \bullet (\alpha \odot \neg \beta)$
\hspace{2 cm} {\footnotesize by Ax P\L 5}

\item[(3)]
$\vdash ((\gamma\bullet \alpha) \odot \neg(\gamma \bullet \beta))
\rightarrow (\alpha \odot \neg \beta)$ \hspace{2.5 cm}{\footnotesize
by 1,2 \L 2}

\item[(4)]
$\vdash (((\gamma\bullet \alpha) \odot \neg(\gamma \bullet
\beta))\rightarrow (\alpha \odot \neg \beta))   \rightarrow (\neg
(\alpha \odot \neg \beta) \rightarrow \neg ((\gamma\bullet \alpha)
\odot \neg(\gamma \bullet \beta)))$ \noindent

\hspace{8.3 cm} {\footnotesize by \L 3 }

\item[(5)]
$\vdash \neg (\alpha \odot \neg \beta) \rightarrow \neg
((\gamma\bullet \alpha) \odot \neg(\gamma \bullet \beta))$
\hspace{1.5 cm} {\footnotesize by MP 3,4}

\item[(6)]
$\vdash (\alpha \rightarrow \beta)\rightarrow  \neg(\alpha \odot
\neg \beta)$ \hspace{3.7 cm} {\footnotesize by \cite[Lemma
3.1.1]{CDM}}

\item[(7)]
$\vdash (\alpha \rightarrow \beta)\rightarrow \neg ((\gamma\bullet
\alpha) \odot \neg(\gamma \bullet \beta))$ \hspace{2.1
cm}{\footnotesize by 5,6 , \L 2}

\item[(8)]
$\vdash \neg ((\gamma\bullet \alpha) \odot \neg(\gamma \bullet
\beta)) \rightarrow ((\gamma\bullet \alpha) \rightarrow (\gamma
\bullet \beta))$ \hspace{0.2 cm} {\footnotesize by \cite[Lemma
3.1.1]{CDM}}

\item[(9)]
$\vdash (\alpha \rightarrow \beta)\rightarrow ((\gamma\bullet
\alpha) \rightarrow (\gamma \bullet \beta))$ \hspace{2.4 cm}
{\footnotesize by MP 3,4}

\end{enumerate}
\qed
\end{proof}

From the above follows that $P\L$ is an admissible expansion of
$\L$. Consequently we can apply Theorem \ref{COMP} to obtain a Pavelka style
completeness for product \L ukasiewicz logic.
The Pavelka style completeness for product \L ukasiewicz logic was obtained in a
different form in {\rm \cite{HC}}.

\subsection{\bf Division \L ukasiewicz logic}

\begin{definition}

{\rm A divisible MV-algebra {\rm \cite{GER}} (for short:
DMV-algebra) is an algebra $ \langle A,  (\delta_n)_{n\in N}
\rangle$ with $\delta_n$ unary operation satisfying the following
for each $x\in A$, $n \in N$:

\begin{enumerate}
\item[1]
$A$ is an MV-algebra

\item[2]
$n\delta_nx= x$

\item[3]
$\delta_nx \odot (n-1)\delta_nx = 0$

\end{enumerate}

}
\end{definition}

We denote by ${\cal DMV}$ the variety of DMV-algebras. If $A$ is a
DMV-algebra then {\rm  \cite[Proposition 5.1.7]{GER}} $Con_{\cal
DMV}(A) = Con_{\cal MV}(A)$, resulting  ${\cal DMV}$ a compatible
expansion of ${\cal MV}$. An important example of DMV-algebra is
$[0,1]_{MV}$ equiped with the unary operations $\delta_n x = x/n$,
the division by $n$. This algebra is denoted by $[0,1]_{DMV}$
containing $Q\cap [0,1]_{DMV}$ as sub algebra. It is clear that
$[0,1]_{DMV}$ is hereditarily simple algebra. By {\rm
\cite[Proposition 5.1.3]{GER}} $(\delta_n = x/n)_{n\in N}$ are the
unique operations making $[0,1]_{MV}$ and $Q\cap [0,1]_{MV}$ into a
DMV-algebra resulting canonical expansions.

\begin{definition}
{\rm We define the division \L ukasiewicz propositional calculus
$D\L$ adding into the language of $\L$ unary connectives
$(\delta_n)_{n\in N} $ plus definig inductively $2\alpha = \neg
\alpha \rightarrow \alpha$, $(k+1)\alpha = \neg \alpha \rightarrow k
\alpha$ and considering the following formulas as axioms

\begin{enumerate}

\item[D\L 0]
\L ukasiewicz axioms

\item[D\L 1]
$ \alpha \Longleftrightarrow k(\delta_k \alpha)$

\item[D\L 2]
$ (\alpha \rightarrow k \beta) \rightarrow (\delta_k \alpha
\rightarrow \beta)$

\end{enumerate}

\noindent The unique deduction rule is the {\it Modus Ponens} }
\end{definition}

\begin{prop} \label{ALGDMV}
{\rm  \cite[Proposition 5.2.1]{GER}}$D\L$ calculus is algebraizable
in ${\cal DMV}$.

\qed
\end{prop}

From the above follows that $D\L$ is an admissible expansion of
$\L$. Consequently we can apply Theorem \ref{COMP} to obtain a Pavelka style
completeness for this logic.

\subsection{\bf \L ukasiewicz logics with the negation fixpoint}
Let ${\cal A}$ be anyone of the following varieties: ${\cal MV}$, ${\cal PMV}$
or ${\cal DMV}$ and $L({\cal A})$ be the corresponding propositional calculus.
Let ${\cal A}_{1/2}$ be the variety built from ${\cal A}$, adjoining the new
constant symbol $k$ and the equation $\neg k = k$. It is clear that if $A$ is a
${\cal A}_{k}$-algebra then, $Con_{{\cal A}_{k}}(A) = Con_{\cal MV}(A)$,
resulting  ${\cal A}_{k}$ a compatible expansion of ${\cal MV}$. $[0,1]_{{\cal
A}}$ equipped with the constant $k = 1/2$ is a ${\cal A}_{k}$-algebra noted
$[0,1]_{{\cal A}_{k}}$. It contains $Q_{[0,1]}$ equipped with the ${{\cal
A}_{k}}$-operations  as sub ${\cal A}_{k}$-algebra. It is clear that 
$[0,1]_{{\cal A}_{k}}$ is hereditarily simple. By {\rm  \cite[Lemma 2.10]{UHO}},
$k=1/2$ is the unique fixpoint of the negation of $[0,1]_{{\cal A}}$. Thus
$[0,1]_{{\cal A}_{k}}$ is the ${\cal A}_{k}$-canonical expansion of
$[0,1]_{MV}$.
We define the propositional calculus $L({\cal A}_{k})$, adding into the language
of $L({\cal A})$ the constant symbol $k$, and consider the following formulas as
axioms:

\begin{enumerate}
\item
$L({\cal A})$-axioms

\item
$\neg k \Longleftrightarrow k$

\end{enumerate}

and modus ponens as unique deduction rule. It is clear that $L({\cal
A}_{k})$ is algebraizable in ${\cal A}_{k}$. From this, follows that
$L({\cal A}_{k})$ is an admissible expansion of $\L$. Consequently
we can apply Theorem \ref{COMP} to obtain a Pavelka style
completeness for these logics (see also {\rm \cite{EGM}}).

{\small \noindent e-mail: hfreytes@gmail.com} 

\end{document}